\documentclass[a4paper,10pt]{amsart}
\usepackage{amssymb}                                                            
\usepackage{eucal}
\usepackage{amsmath}                                                            
\usepackage{amscd}                                                              
\usepackage[dvips]{color}
\usepackage{multicol}
\usepackage[all]{xy}                                                            
\usepackage{graphicx}
\usepackage{color}
\usepackage{colordvi}
\usepackage{xspace}
\usepackage{tikz}                                                               
\usetikzlibrary{arrows,decorations.markings,plotmarks,decorations.markings}
\usetikzlibrary{positioning}	
\usetikzlibrary{snakes} 
\usepackage{tikz-cd}		

\usepackage{mathrsfs}															

\usepackage[active]{srcltx} 

\makeindex                                                                          

\DeclareMathOperator{\Ad}{Ad}

\DeclareMathOperator{\Res}{Res}

\DeclareMathOperator{\SO}{SO}
\DeclareMathOperator{\Spin}{Spin}

\DeclareMathOperator{\Card}{Card}

\DeclareMathOperator{\im}{im}
\begin{document}

\newtheorem{defi}{Definition}[section]
\newtheorem{resu}{Result}[section]
\newtheorem{thm}{Theorem}[section]
\newtheorem{note}{Note}[section]
\newtheorem{lem}{Lemma}[section]
\newtheorem{prop}{Proposition}[section]
\newtheorem{exe}{Exercise}[section]
\newtheorem{cor}{Corollary}[section]

\newcommand{\td}{\mathrm{d}}
\newcommand{\C}{\mathrm{C}}
\newcommand{\e}{\mathrm{e}}
\newcommand{\id}{\mathrm{id}}
\newcommand{\coker}{\mathrm{coker}}
\newcommand{\divi}{\mathrm{Div}}
\newcommand{\pic}{\mathrm{Pic}}

\newcommand{\FF}{\mathscr{F}}
\newcommand{\GG}{\mathscr{G}}
\newcommand{\oo}{\mathcal{O}}

\newcommand{\mf}{\mathfrak{m}}

\newcommand{\re}{\mathbb{R}}
\newcommand{\Z}{\mathbb{Z}}
\newcommand{\p}{\mathbb{P}}
\newcommand{\cpl}{\mathbb{C}}
\newcommand{\ra}{\mathbb{Q}}
\renewcommand{\mod}{\mathrm{mod~}}
\newcommand{\ff}{\mathbb{F}}

\newcommand{\pde}[1][x]{\frac{\partial}{\partial#1}}							
\newcommand{\urt}{\sqrt{-1}}															

\newcommand{\tr}{\mathrm{tr}}
\newcommand{\sig}{\mathrm{Sig}}
\newcommand{\fa}{\mathfrak{a}}
\newcommand{\fg}{\mathfrak{g}}
\newcommand{\fh}{\mathfrak{h}}
\newcommand{\fk}{\mathfrak{k}}
\newcommand{\fs}{\mathfrak{s}}
\newcommand{\ft}{\mathfrak{t}}
\newcommand{\fu}{\mathfrak{u}}
\newcommand{\fp}{\mathfrak{p}}
\title[Twisted Characters and Signature of hermitian forms of a Representation]{An Application of Twisted Character on Signature of a finite-dimensional Representation of a Real reductive Lie group}

\author{Chengyu Du}
\maketitle

\noindent Abstract: Let $G_\cpl$ be a connected complex reductive Lie group, and $G$ be a real form. Let $(\pi,V)$ be a finite-dimensional irreducible representation of $G$. Assume $(\pi,V)$ admits a $G$ invariant hermitian form. {In \cite{Signature-of-a-rep-of-reductive}, an analog of the Weyl dimension formula that instead computes the signature of the invariant form for finite dimensional representations of complex Lie groups is given. The goal of this paper is to give a short alternate proof as an application of a more general twisted character theory of \cite{Dirac-Index-and-Twisted-Characters}.

\section{Introduction}

Let $G_\cpl$ be a complex connected reductive algebraic group. Let $G$ be the group of real points, which is a real reductive Lie group. The subgroup $G$ can also be understood as the fixed points of an anti-holomorphic involution $\sigma$ of $G_\cpl$. Let $(\pi,X)$ be a finite-dimensional irreducible representation of $G_\cpl$ over $\cpl$. By restricting to $G$, it is also an irreducible representation of $G$. When the highest weight of $(\pi,X)$ satisfies the condition stated in Corollary \ref{equivalent condition for herm exist}, $(\pi,X)$ admits a $G$-invariant hermitian form
\begin{equation}\label{g-invariant herm}
\langle\cdot,\cdot\rangle:X\otimes X\to\cpl,\quad\langle\pi(g)v,\pi(g)u\rangle=\langle v,u\rangle,~\forall v,u\in X,~g\in G.
\end{equation}
Since $X$ is irreducible, this $G$-invariant hermitian form is non-degenerate, and unique up to a nonzero real scalar. The form has a signature $(p(X),q(X))$, in short $(p,q)$ if there is no ambiguity. The signature is unique up to switching the sides of $p$ and $q$ if the form is changed by a negative scalar. Thus the absolute value of the difference of the pair $(p,q)$ is well-defined
\begin{equation}
\sig(X):=|p-q|.
\end{equation}
Because $\dim X$ can be calculated by Weyl's dimension formula, to understand $\sig(X)$ is equivalent to understand the signature of $X$.\\
\indent
A known fact is that $\fg_\cpl$ always has a compact real form \cite{CpctReFcplssLa}. A compact real form $\tau$ of $\fg_\cpl$ is a conjugate linear involution of $\fg_\cpl$. Let $\fu$ be the fixed points of $\tau$. Then the corresponding Lie subgroup $U$ of $G_\cpl$ is compact.
Let $(\pi_U,X_U)$ be a finite-dimensional representation of $U$, then it must be unitary. In fact, given any non-degenerate hermitian form $\langle\cdot,\cdot\rangle$ of $(\pi,X)$, we have a positive definite non-degenerate hermitian form
$$\langle u,v\rangle_U:=\frac{1}{|U|}\int_U\langle\pi(g)v,\pi(g)u\rangle dg.$$
Here the integral is taken over the Haar measure. It is $U$-invariant, hermitian, positive definite, and non-degenerate by its definition.\\
\indent
By take differentials, all fintie-dimensional Lie group representations are also Lie algebra representations and hence modules of the universal enveloping algebra. Let $\fg_\cpl$ be the Lie algebra of $G_\cpl$, and $U(\fg_\cpl)$ be the universal enveloping algebra of $\fg_\cpl$. Then $\sigma$ and $\tau$ induce anti-involutions on $U(\fg_\cpl)$. The fixed point of $\sigma$ in $\fg_\cpl$ is $\fg$, which is also the Lie algebra (over $\re$) of $G$. Fix a Cartan subalgebra $\fh_\cpl$ which is stable under $\sigma$. After conjugating of $G_\cpl$, we may assume $\fh_\cpl$ is also stable under $\tau$. The composition of $\tau$ and $\sigma$ is a Cartan involution $\theta=\tau\sigma=\sigma\tau$ who stabilizes $\fh_\cpl$. Let $\fk_\cpl$ be $+1$-eigenspace of $\theta$ and $\fs_\cpl$ the $-1$-eigenspace of $\theta$. Then $\fg_\cpl=\fk_\cpl\oplus\fs_\cpl$. Let $\fk=\fk_\cpl\cap\fg$ and $\fs=\fs_\cpl\cap\fg$. Let $K$ be the Lie subgroup of $G$ whose Lie algebra is $\fk$. There exists a $\cpl$-bilinear form $\langle\cdot,\cdot\rangle$ on $\fg_\cpl$, whose restriction to $\fk$ is negative definite and whose restriction to $\fs$ is positive definite. For example, this form can be the trace form up to scalar.\\
\indent Now fix an irreducible $G_\cpl$ representation $(\pi,X)$. By restriction, it is a $G$ and $U$ representation. Suppose $(\pi,X)$ admits a $G$-invariant hermitian form $\langle\cdot,\cdot\rangle_G$. By taking differentials, $(\pi,X)$ is a $\fg$-representation and this form is $\sigma$-invariant in the following sense,
$$\langle \pi(Y)v,u\rangle_\sigma=-\langle v,\pi(\sigma(Y)u\rangle_\sigma,~\forall~Y\in \fg.$$
Notice that $\langle\cdot,\cdot\rangle_\sigma$ and $\langle\cdot,\cdot\rangle_G$ are the same hermitian form. We use the subscripts $\sigma$ or $G$ depending on the context. Therefor, the signature of this form is preserved.

The formula of the signature of the $G$-invariant hermitian form is given in Theorem \ref{theorem of the conclution}. We will make use of Dirac index theory to provide an alternative proof of Theorem \ref{theorem of the conclution}. First, we define an action of $\theta$ on $X$ and convert the question to calculating the trace of $\theta$. Second, consider the Clifford algebra of $(\fs,\langle\cdot,\cdot\rangle)$ and the universal enveloping algebra $U(\fg)$ of $\fk$. The Dirac operator  $D$ is defined in $U(\fg)\otimes C(\fs)$. We convert the question into calculating the twisted Dirac index $I(X)$.\\

\section{Preliminaries}

In this section, we briefly state some facts about Lie groups, Lie algebras, Dirac operators and the twisted characters. We also state the main theorem that we will prove in an alternative way.

\subsection{Cartan decomposition and root systems}
Let $G_\cpl$ be a complex connected reductive algebraic group, $G$ be the group of real points. The group $G$ can also be considered as a $G_\cpl$ subgroup, which is the fixed points of an anti-holomorphic involution $\sigma$ of $G_\cpl$. Let $\tau$ be an anti-holomorphic involution of $G_\cpl$ such that the fixed points of $\tau$ is a compact Lie group $U$, i.e. $\tau$ is a Lie algebra involution is a compact real form of $\fg_\cpl$. After compositing with a conjugation of $G_\cpl$, the two involutions $\tau$ and $\sigma$ are commutative. Let $\theta=\tau\sigma=\sigma\tau$. Then $\theta$ is a Cartan involution of $G_\cpl$.\\
\indent
Fix a maximal torus $H_\cpl$ which is stable under $\sigma$ and let $H=G\cap H_\cpl$. Let $\fg_\cpl,\fh_\cpl,\fg,\fh$ be the Lie algebra of $G_\cpl,H_\cpl,G,H$ respectively. Then $\fh_\cpl$ (resp. $\fh$) is a Cartan subalgebra of $\fg_\cpl$ (resp. $\fg$). By taking differentials, the involutions of groups induce involutions of Lie algebras. We use the same notation for typing convenience. Because $\theta$, $\tau$ and $\sigma$ commute, $H_\cpl$ is also stable under $\tau$ and $\theta$. Let $H=G\cap H_\cpl$. For the same reason, $H$ and $G$ are stable under $\tau$ and $\theta$. Differentiation preserves those properties, thus we have the Cartan decomposition by the eigenvalue $\pm1$ of the Cartan involution $\theta$,
$$\fg_\cpl=\fk_\cpl\oplus\fs_\cpl,\qquad \fg=\fk\oplus\fs.$$
Here $\theta$ fixes $\fk_\cpl$ and $\fk$, and acts as $-1$ on $\fs_\cpl$ and $\fs$. The corresponding Cartan decomposition of $G$ is $G=K\exp \fs$, where $K$ is the Lie group of $\fk$. The corresponding decomposition of the Cartan subalgebra is
$$\fh_\cpl=\ft_\cpl\oplus\fa_\cpl,\qquad \fh=\ft\oplus\fa,$$
where $\ft\subset\fk$, $\ft_\cpl\subset\fk_\cpl$, $\fa\subset\fs$ and $\fa_\cpl\subset\fs_\cpl$.\\
\indent
We pick $H_\cpl$ to be maximal compact, which means $\dim \ft$ is maximal among all Cartan subalgebras of $\fg$.\\
\indent For $\fg_\cpl$ we have a root system $\Delta(\fg_\cpl,\fh_\cpl)\subset\fh_\cpl^\ast$ (the linear dual space). Since $H$ is maximally compact, $\ft_\cpl$ is a Cartan subalgebra of $\fk_\cpl$, then we also have a root system $\Delta(\fk_\cpl,\ft_\cpl)\subset\ft^\ast_\cpl$. Indeed these two root systems are in different dual spaces. By assigning value 0 on $\mathfrak a_\cpl$ for each function $f\in\ft_\cpl^\ast$, all $f$ can be extended as an element in $\fh^\ast_\cpl$. In this way, we define an embedding, $\ft^\ast_\cpl\hookrightarrow\fh^\ast_\cpl$, and then put $\Delta(\fk_\cpl,\ft_\cpl)\subset\fh^\ast_\cpl$. Besides these two root system, there is the restricted root system $\Delta_R(\fg_\cpl,\fh_\cpl)$ which is obtained by taking restrictions on roots $\Delta(\fg_\cpl,\fh_\cpl)$. In detail, we construct the map $\Res(\alpha):=\alpha|_{\ft_\cpl}\in\ft^\ast_\cpl\subset\fh^\ast_\cpl$, and $\Delta_R(\fg_\cpl,\fh_\cpl)$ is the image of $\Delta(\fg_\cpl,\fh_\cpl)$ through this map. Note that $\Delta_R(\fg_\cpl,\fh_\cpl)$ may not be identical to $\Delta(\fk_\cpl,\ft_\cpl)$, not even necessarily a subset. This restriction map can be realized by the Cartan involution by the following method. Let $\alpha$ be a root of $\Delta(\fg_\cpl,\fh_\cpl)$, $X_\alpha$ be the corresponding root vector. Then $\frac{1}{2}(X_\alpha+\theta X_\alpha)$, if nonzero, is a $\Res(\alpha)$-weight vector living in $\fk_\cpl$. The collection of these vectors exhausts all weight vectors of $\Delta(\fk_\cpl,\ft_\cpl)$. Meanwhile $\frac{1}{2}(X_\alpha-\theta X_\alpha)$, if nonzero, is also a $\Res(\alpha)$-weight vector but living in $\fs_\cpl$. The collection of both the two types of vectors, $\frac{1}{2}(X_\alpha\pm\theta X_\alpha)$, exhaust all weight vectors of $\Delta_R(\fg_\cpl,\fh_\cpl)$. Let $\Delta^+(\fg_\cpl,\fh_\cpl)$ be a $\theta$-stable set of positve roots. We say a restricted root $\Res(\alpha)$ is positive if $\alpha$ is positve. Denote by $W$ the Weyl group of $\Delta(\fg_\cpl,\fh_\cpl)$, denote by $\rho_G$ the half sum of positive roots from $\Delta(\fg_\cpl,\fh_\cpl)$, and denote by $\rho_K$
the half-sum of positive roots from $\Delta(\fk_\cpl,\ft_\cpl)$.
\subsection{The Clifford algebra and spin modules}
\indent Let $V$ be an $n$-dimensional real vector space with an inner product $\langle\cdot,\cdot\rangle$. The Clifford algebra $C(V)$ over $V$ is an associative real algebra with unity, together with a canonical map $i:V\to C(V)$, such that the following universal property holds. Let $A$ be any associative real algebra with unity and let $\phi:V\to A$ be a linear map such that
$$(\phi(v))^2=-\langle v,v\rangle,~\forall v\in V$$
in $A$. One can also consider $C(V)$ as the quotient algebra of $T(V)$ where the kernel is generated by all of the following terms
$$v\otimes w+w\otimes v+2\langle v,w\rangle, ~v,w\in V$$
\indent Fix a non-degenerate, $\mathrm{Ad}(G)$-invariant, and $\theta$-invariant symmetric bilinear form $B$ which is positive definite on $\fs$ and negative definite on $\fk$. If $G$ is semisimple, then the Killing form meets all requirements; If $G$ is a reductive Lie group of matrices, then the trace form would work. Let $C(\fs)$ be the Clifford algebra on $\fs$ with respect to $B$. In the group of units of $C(\fs)$, there lives the spin group $\Spin(\fs)$, who is a double cover of $\mathrm{SO}(\fs)$. We call the covering map $\varphi:\Spin(\fs)\to\SO(\fs)$. Meanwhile, since $B$ is $\mathrm{Ad}(G)$-invariant, there is a canonical group homomorphism $\mathrm{Ad}: K\to \mathrm{SO}(\fs)$. Following the diagram, we have the pullback $\widetilde{K}$, the spin double cover of $K$.

\[\begin{tikzcd}
\widetilde K \arrow[r]\arrow[d]&
\Spin(\fs)\arrow[d,"\varphi"]&\\
K\arrow[r,"\phi"]&
\SO(\fs)
\end{tikzcd}\]

One of the objects we need in this paperis the complex spin module. In \cite[Chapter 2]{DiracOperatorsHuangJingsongBook}, the book shows the $C(\fs)$-module structure of the spin module $S$, and hence $S$ is a representation of the spin group $\Spin(\fs)$ by restricting the action of $C(\fs)$ to $\Spin(\fs)$. The differential of the action of $\Spin(\fs)$ corresponds to the action of $\mathfrak{so}(\fs)\subset C(\fs)$. Here we briefly show how the complex spin module looks like as a vector space. Notice that $\fg_\cpl$ can be decomposed as a direct sum of $\Delta(\fk_\cpl,\ft_\cpl)$-weight subspaces. We pick $\{u_i\}$ to be nonzero weight vectors of the positive weights $\{\beta_i\}$, and $\{u_i^\ast\}$ for negative weights $\{-\beta_i\}$. By $\Ad(G)$-invariance of $B$ and non-degeneracy of $\fs_\cpl$, it easily follows that $(u_i,u_i^\ast)$ are non-degenerately paired. All zero weights live in $\fa_\cpl$ (because we chose $H_\cpl$ to be maximally compact). When $\dim\fs_\cpl$ is even, $\dim\fa_\cpl$ is also even, and $\fa_\cpl$ is a direct sum of two isotropic subspaces dual to each other. Let $\{v_j\}$ be a basis of one of the isotropic subspaces. Then $\{u_i,v_j\}$ spans an isotropic subspace $L$ of $\fs_\cpl$. The spin module $S$ as a vector space is $\bigwedge L$, whose dimension is $2^{\dim\fs/2}$. Let $(I,J)$ be a pair of index sets, $u_I:=\bigwedge_{i\in I}u_i$ and $v_J:=\bigwedge_{j\in J}v_j$. Then $u_I\wedge v_J$ is a $\Delta(\fk_\cpl,\ft_\cpl)$-weight vector via $d\phi$, and the weight depends on $I$ but is independent of $J$. So a fixed weight has multiplicity $2^{\dim \fa/2}$. When $\dim_\cpl\fs$ is odd, $\dim\fa$ is also odd. We could make similar construction as the even dimensional case. In this case, the dimension of the spin module is $2^{(\dim\fs-1)/2}$, and each weight corresponding to a fixed $I$ has multiplicity $2^{(\dim\fa-1)/2}$.

\subsection{The Dirac operator and the main theorem }
Let $(\pi,V)$ be a representation of $G_\cpl$. It's then a $(\fg,K)$-module. Make a tensor product of the $(\fg,K)$-module $V$ and the $(C(\fs),\Spin(\fs))$-module together, then we have a $U(\fg)\otimes C(\fs)$-module $V\otimes S$. Noitce $\widetilde{K}$ is a pullback of $K$ and $\Spin(\fs)$, we could write $\widetilde{K}=\{(k,g)|\phi(k)=\varphi(g),~k\in K,~s\in\Spin(\fs)\}$. Then $V\otimes S$ is a representation of $\widetilde{K}$. The actions of $\widetilde K$ is defined in the following way. Given an element $(k,g)\in\widetilde K$, $k$ acts on $V$ via the representation of $G$ and $g$ acts on $S$ via the representation of $\Spin(\fs)$.\\
\indent
\begin{defi}\rm\cite{DiracOperatorsHuangJingsongBook}
	Let $\{Z_i\}$ be an orthonormal basis of $\fs$, the Dirac operator $D$ is an element in $U(\fg)\otimes C(\fs)$, defined as
	$$D=\sum_i Z_i\otimes Z_i.$$
	It's easy to see that $D$ is independent of the choice of the basis $\{Z_i\}$ as long as they are orthonormal.
\end{defi}
Given a $(\fg,K)$-module $V$ and a spin module $S$, the Dirac operator $D$ acts on the module $V\otimes S$. An important property is that the action of $D$ commutes with the action of $\widetilde{K}$ on $V\otimes S$. Hence $\ker D$ and $\im D$ are representations of $\widetilde{K}$.
\begin{defi}\rm\cite{DiracOperatorsHuangJingsongBook}
	Given a $(\fg,K)$-module $V$, and fix a spin module $S$, the Dirac operator acts on $V\otimes S$. The Dirac cohomology $H^D(V)$ is defined as
	\begin{equation}
	H^D(V)=\ker D/(\ker D\cap\im D).
	\end{equation}
\end{defi}
When $\dim V<\infty$ or $V$ is unitary, the operator $D$ satisfies $\ker D^2=\ker D$. An direct corollary is that $\ker D\cap\im D=0$ and hence $H_D(V)=\ker D$.\\
\indent
We also need the Vogan's conjecture (also proved).
\begin{thm}\rm\cite{DiracOperatorsHuangJingsongBook}
	Let $V$ be an irreducible $(\fg,K)$-module. Assume that the Dirac cohomology of $V$ contains a $K$-type $E(\mu)$ of highest weight $\mu\in\ft^\ast\subset\fh^\ast$. Then the infinitesimal character of $V$ is $\mu+\rho_K$.
\end{thm}

A direct corollary of Vogan's conjecture is that $H_D(V)$ can be decomposed into irreducible $\widetilde K$-types $E^{\widetilde K}(\mu)$ such that $\mu$ has the form
$$w(\lambda+\rho_G)-\rho_K,w\in W^1,$$
where $W^1=\{w\theta=\theta w|w(\Delta_R^+(\fg_\cpl,\fh_\cpl))\supset\Delta^+(\fk_\cpl,\ft_\cpl)\}$. The highest weight vectors have the shape $z_{w\lambda}\otimes u_I\wedge v_J,~w\in W^1$. So each irreducible type has multiplicity $2^{\dim\fa/2}$ or $2^{(\dim\fa-1)/2}$. See the detaisl in \cite{Signature-of-a-rep-of-reductive}.\\
\indent
Before we proceed, we need to introduce a new number $\epsilon(w)$. For any element $w\in W^1$, we could write $w\lambda$ in the form
\begin{equation}\label{how w lambda looks like}
w\lambda=\lambda-\sum_\beta (n_\beta\cdot\beta)
\end{equation}
where $\beta$ must be singular (here is means perpendicular to the sum of complex roots) and imaginary. Those numbers $n_\beta$ help computing $\epsilon(w)$. \cite[Corollary 4.3]{Signature-of-a-rep-of-reductive}, we pick up the $\beta$'s whose are noncompact imaginary, and then $\epsilon(w)$ is defined to be $\displaystyle\prod_{\beta~\mathrm{ncpt}}(-1)^{n_\beta}$.\\
\indent Now it's a time to state the main theorem in this paper. Let $(\pi,V)$ be a finite dimensional irreducible representation of $G_\cpl$. Assume it admits an $G$-invariant hermitian form. Let $\lambda$ be the highest weight of $V$. Given any non-degenerate hermitian form, its signature is presented as a pair $(p,q)$. Define
$$\sig(V)=|p-q|.$$
Because $p+q=\dim V$ can be calculated by the Weyl's dimension formula, knowing their different $|p-q|$ can explicitly find $p$ and $q$. The following theorem gives the formula for $\sig(V)$.
\begin{thm}\rm\label{theorem of the conclution}\cite{Signature-of-a-rep-of-reductive}
	Let $V$ be an irreducible finite-dimensional representation of $G$. Assume $V$ has highest weight $\lambda$. Then
	\begin{equation}\label{conclution}
	\sig (V)=\left|\sum_{w\in W^1}(-1)^{\epsilon(w)}\dim E^{\tilde{K}}(w(\lambda+\rho_G)-\rho_K)/2^r\right|.
	\end{equation}
	Here $E^{\tilde K}(\nu)$ is a finite-dimensional representation of $\tilde{K}$ and $\nu$ is its highest weight.
\end{thm}
\indent
In this paper, we will apply the twisted character theory in \cite{Dirac-Index-and-Twisted-Characters} to give an alternate proof of Theorem \ref{theorem of the conclution}. Given a character formula of an irreducible group representation, its value at identity is the dimension. In our approach, we use a character twisted by an automorphism of the group $\widetilde K$. We evaluate this twisted character not only at the identity, but also an element in the extended group to obtain the key equations.

\subsection{A Brief Introduction to the Twisted Character}

This subsection will show what we need from \cite{Dirac-Index-and-Twisted-Characters}. Let $\gamma$ be an automorphism of $(U(\fg_\cpl)\otimes C(\fs_\cpl),\tilde K)$. This means the following:\\
(1) $\gamma$ consists of an automorphism $\gamma_\fg$ of $U(\fg_\cpl)\otimes C(\fs_\cpl)$ and an automorphism $\gamma_K$ of $\tilde K$;\\
(2) $\gamma$ is compatible with the action of $\tilde K$ on $U(\fg_\cpl)\otimes C(\fs_\cpl)$ in the sense that
$$\gamma_\fg((k, g)(u\otimes c))=\gamma_K (k, g)\gamma_\fg(u\otimes c), ~\forall~(k,g)\in \tilde K~\mathrm{and}~u\otimes c\in U(\fg_\cpl)\otimes C(\fs_\cpl);$$
(3) the differential of $\gamma_K$ coincides with the restriction of $\gamma_\fg$ to $\fk_\Delta$. Here $\fk_\Delta$ means the image of $X\mapsto X\otimes 1+1\otimes d\phi(X)$, where $d\phi$ is the differential of the map $\phi: K\to\SO(\fs)$ followed by complexification.

\indent Let $V\otimes S$ be a $(U(\fg_\cpl)\otimes C(\fs_\cpl),\tilde K)$-module as above. We assume that $V\otimes S$ has a compatible action of $\gamma$, i.e. $\gamma$ is an operator on $V\otimes S$ such that
\begin{gather*}
\gamma\circ(u\otimes c)\circ\gamma^{-1}=\gamma_\fg(u\otimes c),\forall~u\otimes c\in U(\fg_\cpl)\otimes C(\fs_\cpl).\\
\gamma\circ k\circ\gamma^{-1}=\gamma_K(k),\forall~ k\in\tilde K.
\end{gather*}
The tools we apply for our proof are the following equations. Let $\tilde K_\gamma$ be the set of $\gamma$-fixed points in $\tilde K$. We use $+$ and $-$ on the upper right corner to denote the $+1$ and $-1$ eigenspace of $\gamma$. Then the twisted character is defined as 
\begin{equation}\label{character on D-coho}
\chi_\gamma^V(k)=\tr(k;H_D(V)^+)-\tr(k;H_D(V)^-),\forall~k\in \tilde K_\gamma
\end{equation}
It could also be equivalently considered as a virtual $\tilde{K}_\gamma$-module
\begin{equation}
I_\gamma(V)=H_D(V)^+-H_D^-(V).
\end{equation}
Later in \cite[Proposition, 1.8]{Dirac-Index-and-Twisted-Characters}, we have
\begin{equation}
I_\gamma(V)=(V\otimes S)^+-(V\otimes S)^-.
\end{equation}
As character function, it is
\begin{equation}\label{character on X and S}
\chi_\gamma^V(k)=\tr(k;(V\otimes S)^+)-\tr(k;(V\otimes S)^-),~\forall k\in\tilde K_\gamma
\end{equation}
Combine (\ref{character on D-coho}) and (\ref{character on X and S}), we will have
\begin{equation}\label{frequently used formula}
\tr(k;H_D(V)^+)-\tr(k;H_D(V)^-)=\tr(k;(V\otimes S)^+)-\tr(k;(V\otimes S)^-).
\end{equation}
This equation will be the one that builds up our proof in the next section.

\section{Our Approch}
In the following discussions, let $G_\cpl$ be a connected complex reductive Lie group, $G$ be a real form, and $(\pi,V)$ be a finite-dimensional representation of $G_\cpl$. All notations in previous sections are valid in this section. If no ambiguities exist, then we also regard $(\pi,V)$ as a representation of $G$, or the corresponding representation of Lie algebras $\fg_\cpl$ and $\fg$. And sometimes for short we write $X.v$ to indicate $\pi(X)v$ where $X\in\fg_\cpl,~v\in V$.\\

\subsection{The action of $\theta$ on the underlining vector space $V$}
\begin{defi}\rm
	Given $(\pi,V)$ and the compact real form $\tau$ as above, $(\pi^\theta,V)$ is defined to be a representation of $G_\cpl$ equipped with exactly the same vector space structure, but the group homomorphism is defined as
	\begin{equation}
		\pi^\theta(g)=\pi(\theta g),~\forall g\in G_\cpl
	\end{equation}
	By taking differential, $(\pi^\theta,V)$ is also a representation of $\fg_\cpl$. It satisfies
	$$\pi^\theta(X)=\pi(\theta X),~\forall X\in\fg_\cpl.$$
	One can easily check that this is well-defined. 
\end{defi}

\begin{lem}\rm
	Given $(\pi,V)$ as above, there exists a positive definite hermitian form $\langle\cdot,\cdot\rangle_u$ satisfying
	\begin{equation}\label{tao inv}
		\langle X.v,w\rangle_u=\langle v,-\tau(X).w\rangle_u, \forall X\in\fg_\cpl,~v,w\in V.
	\end{equation}
	We say that a form satisfying Eq. \ref{tao inv} is $\tau$-invariant.
	\begin{proof}
		Let $U$ be the compact form of $G$. Since $\dim V<\infty$, there exists a $U$-invariant positively definite nondegenerate hermitian form $\langle\cdot,\cdot\rangle_u$. By doing differential we obtain a representation of the Lie algebra $\mathfrak u$, which admits an invariant hermitian form
		$$\langle X.v,w\rangle_u=\langle v,-X.w\rangle_u, \forall X\in\mathfrak u,~v,w\in V.$$
		Then by complexification of $\mathfrak u$, we obtain the $\tau$-invariant form on $V$.
	\end{proof}
\end{lem}

\begin{lem}\rm\label{sigma inv}
	Given $(\pi,V)$ as stated. Let $\langle\cdot,\cdot\rangle$ be a hermitian form on $V$. Then $\langle\cdot,\cdot\rangle$ is $\sigma$-invariant if only if it is $G$-invariant. By $\sigma$-invariance we mean Eq. \ref{tao inv} with $\tau$ replaced by $\sigma$.
	\begin{proof}
		Suppose $\langle\cdot,\cdot\rangle$ is $G$-invariant. It is straightforward via differentiation to see that the form it is $\sigma$-invariant.\\
	\indent Suppose $\langle\cdot,\cdot\rangle$ is $\sigma$-invariant. Notice $G_\cpl$ is connected, and $\dim V<\infty$. The $\sigma$-invariance can be lifted to a group property, which is the $G$-invariance.
	\end{proof}
\end{lem}

\begin{prop}\rm\label{eqiuvalent condition of existence of H form and theta action}
	Given $(\pi,V)$ as above, the following two statement are equivalent.\\
	(a) $V$ admits a $G$-invariant nondegenerate hermitian form $\langle\cdot,\cdot\rangle$.
	$$\langle X.v,w\rangle=\langle v,-\sigma(X).w\rangle, \forall X\in\fg_\cpl,~v,w\in V.$$
	(b) The two representation $(\pi,V)$ and $(\pi^\theta,V)$ are isomorphic as representations of $G$. That is, there exists a map $T:V\to V$, and it satisfies $T\circ \pi(g)=\pi(\theta g)\circ T$.
	\begin{proof}
		Suppose (b) holds. Define $\langle v,w\rangle:=\langle Tv,w\rangle_u$. By definition it's nondegenerate and hermitian. Then
		$$\langle X.v,w\rangle=\langle T(X.v),w\rangle_u
		=\langle (\theta X).(Tv),w\rangle_u
		=\langle Tv,-(\tau(\theta X)).w\rangle_u
		=\langle Tv,-(\sigma X).w\rangle_u
		=\langle v,-(\sigma X).w\rangle.$$
		So this newly defined hermitian form is $\sigma$-invariant. Lemma \ref{sigma inv} lifts the hermitian form to a $G$-invariant one.\\
	\indent	Suppose (a) holds. Let $v\in V$, define $Tv\in V$ to be the unique element such that $\langle v,w\rangle=\langle Tv,w\rangle_u$. $T$ is a linear isomorphism since both forms are nondegenerate. It's easy to check that $T$ is an isomorphism of the Lie algebra $\fg_\cpl$ representations, i.e. $T\circ \pi(X)=\pi(\theta X)\circ T$. The arithmetic work is similar as the computation above. Then by the connectedness of $G_\cpl$ and the fact $\dim V<\infty$, this commuting property can be lifted to the Lie group $G$.
	\end{proof}
\end{prop}

\begin{cor}\rm
	If the Cartan involution $\theta$ is inner, i.e. $\theta=\Ad(k_0)$ for some $k_0\in K$. Then any representation $(\pi,V)$ of $G$ admits a $G$-invariant hermitian form.
	\begin{proof}
		Indeed $\pi(k_0)$ exactly plays the role of $T$ in the proof of Proposition \ref{eqiuvalent condition of existence of H form and theta action}.
	\end{proof}
\end{cor}

\begin{cor}\rm\label{equivalent condition for herm exist}
	Consider the decomposition $\fh^\ast_\cpl=\ft^\ast_\cpl\oplus\fa^\ast_\cpl$. Let $V$ be a finite dimensional irreducible representation of $G$ with highest $\lambda=(\lambda_t,\lambda_a)$, where $\lambda_t\in\ft^\ast_\cpl$, $\lambda_a\in\fa_\cpl^\ast$. Then $V$ admits a $G$-invariant hermitian form if and only if $\theta\lambda=(\lambda_t,-\lambda_a)$ is $W$-conjugate to $\lambda$. In particular, if the Cartan subgroup $H$ is maximally compact, and $K$ has the same rank as $G$. Then $V$ admits a $G$-invariant hermitian form.
	\begin{proof}
		It's sufficient to consider the highest weight of $(\pi^\theta,V)$. Let $v$ be a highest weight vector of $(\pi,V)$, then $v$ is of weight $\theta\lambda$. It follows that $\theta\lambda$ is an extremal weight of $(\pi^\theta,V)$. The equivalent condition given by Proposition \ref{equivalent condition for herm exist} is that $(\pi,V)$ and $(\pi^\theta,V)$ are isomorphic as $G$-representation, which is equivalent to say that they have the highest weight. It is then equivalent to say $\theta\lambda$ as an extremal weight is conjugate to $W$.\\
\indent
		When $H$ is maximally compact and $K$ has the same rank as $G$ does, $\fa_\cpl^\ast=0$. Then $\theta\lambda=\lambda$.
	\end{proof}
\end{cor}

\begin{lem}\rm
	If the isomorphism $T$ exists, then $T^2$ acts as a scalar on $V$.
	\begin{proof}
		By the definition of $T$, it satisfies $T\circ \pi(g)=\pi(\theta g)\circ T$. A simple subtitution can prove $T\circ \pi(\theta g)=\pi(g)\circ T$. Then
		$$T^2\pi(g)=T\pi(\theta g)T=\pi(g)T^2.$$
		Then by Schur's lemma, $T^2$ acts as scalar on $V$.
	\end{proof}
\end{lem}

Without loss of generality, we assume $T^2=\id$, so $T$ is an involution. With such an assumption, we define $\theta$ action on $V$, by setting $\theta v=Tv$. $\theta$ has two eigenvalues $\pm1$. Let $x$ be a $+1$ eigenvector, and $y$ a $-1$ eigenvector. Then
\begin{gather*}
	\langle x,x\rangle=\langle\theta x,x\rangle_u=\langle x,x\rangle_u>0.\\
	\langle y,y\rangle=\langle\theta y,y\rangle_u=\langle -y,y\rangle_u<0.
\end{gather*}
A direct corollary is that the eigenvalues of $\theta$ tells the subspaces where $\langle\cdot,\cdot\rangle$ is positive or negative definite. Hence we have
\begin{equation}
	\sig (V)=|\tr(\theta;V)|.
\end{equation}

\subsection{The action of $\theta$ on the spin module $S$, $U(\fg_\cpl)$ and $C(\fs_\cpl)$}
To define $\theta$-action on $S$, it's sufficient to define $\theta$-action on the basis $\{u_I\wedge v_J\}$.
$$\theta(u_I\wedge v_J)=(-1)^{\Card I+\Card J}u_I\wedge v_J.$$
Let $U(\fg_\cpl)$ be the universal enveloping algebra of $\fg_\cpl$. The universal property directly extends $\theta$-action on $\fg_\cpl$ to $U(\fg_\cpl)$. $\theta$ acts on $\fs_\cpl$ as $-1$. Then the universal property of $C(\fs_\cpl)$ also extends $\theta$-action to $C(\fs_\cpl)$. We can regard $C(\fs_\cpl)$ as a $\Z_2$ graded algebra, with all $\fs_\cpl$ elements odd. Then $\theta$ is $1$ on even part, and $-1$ on odd part.

\subsection{For equal rank case $\dim\ft=\dim\fh$} In \cite{Dirac-Index-and-Twisted-Characters}, $\gamma$ is taken to be $\gamma_1=\id$ on $(\fg_\cpl,K)$, and $\gamma_2=\theta$ on $C(\fs_\cpl)$. Then $\tilde K_\gamma=\tilde K$ which means $\gamma$ commutes with $\tilde K$ actions. Combine Eq \ref{character on D-coho} and \ref{character on X and S}, we have
\begin{equation}\label{equal rank case}
	\tr(k;H_D(V)^+)-\tr(k;H_D(V)^-)=\tr(k;(V\otimes S)^+)-\tr(k;(V\otimes S)^-).
\end{equation}
The tricky part is to pick up a suitable $k\in\tilde{K}_\gamma$. Since $\theta$ is an inner automorphism of $\fg$, there exists $k_0\in K$ s.t. $\theta=\Ad k_0$. And $\theta$ is also an inner automorphism of $C(\fs_\cpl)$ because $\dim\fs$ is even for equal rank case. So $\theta\in\Spin(\fs)$. Both of there image in $SO(\fs)$ is $-\id$. Thus $k_0$ has a preimage $\tilde k_0$ in $\tilde K_\gamma$, and its action on $V\otimes S$, $U(\fg_\cpl)$ and $C(\fs_\cpl)$ exactly as $\theta$. For writing convenience, we use $\theta$ as the notation of $\tilde k_0$. The element we want to plug into Eq. (\ref{equal rank case}) is $\theta$.\\

\noindent\emph{Left Hand Side of Equation (\ref{equal rank case}).}\quad Since $\gamma$ commutes with $\tilde K$ action on $V\otimes S$, which contains $H_D(V)^\pm$, both of this subspaces can be decomposed as a direct sum of $\tilde K$-types $E^{\tilde K}(\mu)$. Because $\theta$ commutes with $\tilde K$ action on $V\otimes S$, so on an irreducible $\tilde K$-type $E^{\tilde K}(\mu)$, $\theta$ has only one eigenvalue. This means it's sufficient to check the eigenvalue of $\theta$ on each highest weight.\\
\indent The highest weight vector of $E^{\tilde K}(\mu)=E^{\tilde K}(w(\lambda+\rho_G)-\rho_H)$ has the form $z_{w\lambda}\otimes u_I\wedge u_J$, as discussed previous sections. For equal rank case, $\dim\fa=0$. Thus $J$ is empty. If $z_{w\lambda}\otimes u_I\in H_D(V)^+$, by the construction of $\gamma$, we see $\theta$ has eigenvalue 1 on $u_I$. The trace of $\theta$ is determined by its eigenvalue on $z_{w\lambda}$.
Recall how $w\lambda$ looks like in Eq. (\ref{how w lambda looks like}), $w\lambda=\lambda-\sum_\beta (n_\beta\cdot\beta)$, and $\beta$ are all imaginary. Up to a difference of a nonzero scalar, we have
$$z=\prod_\beta(X_{-\beta}^{n_\beta}).z_\lambda.$$
Let's pick up a particular $X_{-\beta}$ acting on a weight vector $u$ who is also an eigenvector of $\theta$, and see how $X_{-\beta}$ affects the $\theta$-eigenvalue.
$$\theta(X_{-\beta}.u)=(\theta X_{-\beta}).(\theta u)=\begin{cases}
X_{-\beta}.\theta u&\beta~\text{compact},\\
-X_{-\beta}.\theta u&\beta~\text{noncompact}.
\end{cases}$$
Without loss of generality, we assume $\theta$ acts as $1$ on $z_\lambda$,the highest weight vector of the $(\fg,K)$-module $V$, then by induction one can see that $z_{w\lambda}$ has $\theta$-eigenvalue $\displaystyle\prod_{\beta}(-1)^{n_\beta}$, where $\beta$ are noncompact imaginary roots who show up in Eq. (\ref{how w lambda looks like}). Therefore, this eigenvalue is exactly $(-1)^{\epsilon(w)}$.
We can go over the same discuss on $H_D(V)^-$. The only difference is that for the highest weight vector $z_{w\lambda}\otimes u_I$, $\theta$ acts on $u_I$ as $-1$. For both $H_D(V)^\pm$, the multiplicity of each irreducible type is 1 because $J$ is trivial here.\\
\indent As a conclusion, we have the formula for left hand side of Equation (\ref{equal rank case})
\begin{equation}
	\tr(\theta;H_D(V)^+)-\tr(\theta;H_D(V)^-)=\sum_{w\in W^1}(-1)^{\epsilon(w)}\dim E(w(\lambda+\rho_G)-\rho_K)
\end{equation}

\noindent\emph{Right Hand Side of Equation (\ref{equal rank case}).} Notice $\gamma$ restricted on $S$ is $\theta$. We define $S^{\theta\pm}$ to be the $\pm1$ eigenspace of $\theta$ respectively, then
\begin{align*}
	\tr(\theta;(V\otimes S)^+)-\tr(\theta;(V\otimes S)^-)&=\tr(\theta;V)\cdot\dim S^{\theta+}-\tr(\theta;V)\cdot(-\dim S^{\theta-})\\
	&=\tr(\theta;V)\cdot(\dim S^{\theta+}+\dim S^{\theta-})\\
	&=\tr(\theta;V)\cdot\dim S
\end{align*}
Combine these two equations we have 
\begin{equation}\label{almost conclusion}
	\tr(\theta;V)=\sum_{w\in W^1}(-1)^{\epsilon(w)}\dim E(w(\lambda+\rho_G)-\rho_K)/\dim S.
\end{equation}
For equal rank case, we know $\mathfrak a=0$. Recall $2r=\dim\fs-\dim\fa$, we have $2r=\dim\fs$. Notice that $\dim\fs$ is even, and $\dim S=2^{\dim\fs/2}$. The proof is done.

\subsection{For nonequal rank case, $\dim\ft<\dim\fh$} In \cite{Dirac-Index-and-Twisted-Characters}, $\gamma$ is taken to be $\gamma_1=\theta$ on $(\fg_\cpl,K)$, and $\gamma_2=\id$ on $C(\fs_\cpl)$. With this choice, $\tilde K_\gamma=\tilde K$. For Eq (\ref{character on D-coho}) and (\ref{character on X and S}), pick $k=1$, then we have
\begin{equation}\label{nonequal rank case}
\tr(1;H_D(V)^+)-\tr(1;H_D(V)^-)=\tr(1;(V\otimes S)^+)-\tr(1;(V\otimes S)^-).
\end{equation}
\vskip 6pt
\noindent\emph{Left Hand Side of Equation (\ref{nonequal rank case}).} 
Since we pick $k=1$, the left hand side is a difference of dimensions $\dim H_D(V)^+-\dim H_D(V)^-$ of two $\gamma$-eigenspaces. Notice $\gamma$ commutes with $\tilde{K}$-action, so each $E^{\tilde{K}}$-type is a $\gamma$-eigenspace. Thus it's sufficient to study the $\gamma$-action on those highest weight vector $z_{w\lambda}\otimes u_I\wedge u_J$ for each $E^{\tilde{K}}$-type of highest weight $\mu=w(\lambda+\rho_G)-\rho_K,w\in W^1$. One step further, $\gamma$ acts on $C(\fs_\cpl)$-components as identity, so the eigenvalue is determined by $\theta$-eigenvalue on $z_{w\lambda}$. At this point we could repeat the discussion from the equal-rank case and obtain that $\gamma$ has eigenvalue $(-1)^{\epsilon(w)}$ on $z_{w\lambda}\otimes u_I\wedge u_J$. And hence the left hand side of Equation (\ref{nonequal rank case}) is
$$\sum_{w\in W^1}(-1)^{\epsilon(w)}\dim E^{\tilde{K}}(w(\lambda+\rho_G)-\rho_K)\cdot m$$
As discussed in Section 1, if $\dim \mathfrak a$ is odd, then $m=2^{(\dim\mathfrak a-1)/2}$; if it is even, $m=2^{\dim\mathfrak a/2}$.\\

\noindent\emph{Right Hand Side of Equation (\ref{nonequal rank case}).} Notice $\gamma$ restricted on $V$ is $\theta$. We define $V^{\theta\pm}$ to be the $\pm1$ eigenspace of $\theta$ respectively,
\begin{align*}
	\tr(1;(V\otimes S)^+)-\tr(1;(V\otimes S)^-)&=\dim(V\otimes S)^+-\dim(V\otimes S)^-\\
	&=\dim V^{\theta+}\cdot\dim S+\dim V^{\theta-}\cdot\dim S\\
	&=\tr(\theta;V)\cdot\dim S
\end{align*}
Now combine what we have computed of the left hand and right hand side of Equation (\ref{nonequal rank case}), we have
$$\tr(\theta;V)=\sum_{w\in W^1}(-1)^{\epsilon(w)}\dim E^{\tilde{K}}(w(\lambda+\rho_G)-\rho_K)\cdot m/\dim S.$$
By the construction of spin module,
$$\dim S=\begin{cases}
2^{(\dim\fs-1)/2}& \dim \fs~\text{odd}\\
2^{\dim\fs/2}& \dim \fs~\text{even}
\end{cases}$$
An important fact is that $\dim\fs$ and $\dim\mathfrak a$ have the same parity. So for both cases, $m/\dim S=2^r$. Again we meet Eq. (\ref{almost conclusion}), and the proof is done. 

\subsection{An Example of $SU(3,1)$}

Let $G$ be the group $SU(3,1)$. We compute the signature of it's adjoint representation on the Lie algebra $V=\mathfrak{sl}(4,\cpl)$. The definition of $SU(3,1)$ naturally gives the $G$-invariant hermitian form on $\mathfrak{sl}(4,\cpl)$. It's easy to see that the signature is $(9,6)$, or $(6,9)$ if one reverses the sign of the form. Now we are going to apply Equation (\ref{conclution}) to prove the signature.\\
\indent
The maximal compact subgroup $K$ is isomorphic to $SU(3)\times SU(1)$. The root system of $G$ has three simple roots, two of which are compact and the last one is noncompact. We denote the two compact simple roots by $\alpha,\beta$, the noncompact simple root by $\gamma$. Then $\dim\fk=10$, $\dim\fs=6$ and $\dim\fa=0$. So $2r=\dim\fs-\dim\fa=6$.
The highest weight of the adjoint representation is $\alpha+\beta+\gamma$. The sum of all positve roots is $2\rho_G=3\alpha+4\beta+3\gamma$. The sum of compact positive roots is $2\rho_K=2\alpha+2\beta$. Let $\sigma_\delta$ denote the reflection with respect to a root $\delta$. The set $W^1=\{\id,\sigma_\gamma,\sigma_\gamma\sigma_\beta,\sigma_\gamma\sigma_\beta\sigma_\alpha\}$. 
Apply Equation (\ref{conclution}), we have
\begin{align*}
\sig(V)=&|\dim  E^{\tilde{K}}(\frac{1}{2}(3\alpha+4\beta+5\gamma))+\dim E^{\tilde{K}}(\frac{1}{2}(3\alpha+4\beta+\gamma))\\
&+\dim E^{\tilde{K}}(\frac{1}{2}(-\alpha+2\beta-\gamma))+\dim E^{\tilde{K}}(-\alpha-2\beta-5\gamma)|/2^3\\
=&|3-15-15+3|/8=3.
\end{align*}
Notice that $\dim V=15$. So the signature is $(9,6)$ or $(6,9)$.
\bibliographystyle{plain}
\def \noopsort #1{}


\begin{thebibliography}{1}
	
	\bibitem{Dirac-Index-and-Twisted-Characters}
	Pavle~Pand\u{z}i\'{c} Dan~Barbasch and Peter Trapa.
	\newblock Dirac index and twisted characters.
	\newblock {\em Transactions of American Mathematical Society}, 371:1701--1733,
	2019.
	
	\bibitem{DiracOperatorsHuangJingsongBook}
	Jing-Song Huang and Pavle Pand\u{z}i\'{c}.
	\newblock {\em Dirac Operators in Representation Theory}.
	\newblock Mathematics: Theory and Applications. Birkh\"{a}user Boston, 2006.
	
	\bibitem{Signature-of-a-rep-of-reductive}
	D.~Kalinov, D.~Vogan, Jr., and Christopher Xu.
	\newblock Signatures for finite-dimensional representations of real reductive
	lie groups.
	\newblock {\em arXiv: Representation Theory}, 2018.
	
	\bibitem{CpctReFcplssLa}
	Roger Richardson.
	\newblock {Compact real forms of a complex semi-simple Lie algebra}.
	\newblock {\em Journal of Differential Geometry}, 2(4):411 -- 419, 1968.
	
\end{thebibliography}
\end{document}